# A REMARK ON THE REGULARITY OF PREHOMOGENEOUS VECTOR SPACES

Akihiko Yukie[1]

Oklahoma State University

Let $k$ be a field, $G$ a reductive group over $k$, and $V$ a rational representation of $G$. Note that the reductiveness includes the condition that $G$ is smooth over $k$. Let $e \in G$ be the unit element. Throughout this paper, we assume that there exists $w \in V_k$ such that $U = Gw$ is Zariski open in $V$.

A non-constant polynomial $\Delta(x) \in k[V]$ is called a relative invariant polynomial if there exists a character $\chi$ (which is automatically rational) such that $\Delta(gx) = \chi(g)\Delta(x)$ for all $g \in G, x \in V$. A pair $(G, V)$ as above is sometimes called a prehomogeneous vector space (see [7], [6]), but this terminology is also used with an additional assumption of the existence of a relative invariant polynomial.

Let $G_w$ be the scheme theoretic stabilizer of $w$ in $G$. This is a group scheme and it may not always be reduced. For any scheme $X$ over $k$ and a $k$-algebra $R$, $X_R$ is the set of $R$-rational points of $X$. The purpose of this paper is to prove the following theorem algebraically.

**Theorem** *Suppose $G_w$ is reductive. Then*
*(1) $V \setminus U$ is a hypersurface, and*
*(2) $U_{k^{\mathrm{sep}}}$ is a single $G_{k^{\mathrm{sep}}}$-orbit set-theoretically.*

Note that the assumption in the above theorem includes the smoothness of $G_w$ over the ground filed $k$. The above theorem is well known if ch $k = 0$.

Before proving the above theorem, we recall properties of smooth morphisms. Let $f : X \to Y$ be a morphism of $k$-schemes of finite type, and $\Omega^1_{X/Y}$ the sheaf of relative differentials. Then $f$ is, by definition, smooth of relative dimension $n$ if $f$ is flat and for each point $y \in Y$, the absolute fiber $f^{-1}(y) \times \mathrm{Spec}\,\overline{k(y)}$ ($\overline{k(y)}$ is the algebraic closure) is regular of equidimension $n$. An alternative definition (see Definition (1.1) [1, p. 128]) is that locally $f$ is etale over an affine $n$-space over $Y$. It is easy to see that smoothness is preserved by base change. It satisfies the property of "faithfully flat descent" (see Corollary (5.11) [1, p. 154]). Therefore, if $f : X \to Y$ is a morphism between $k$-schemes of finite type, $f$ is smooth if and only if the induced map $X \times \mathrm{Spec}\,\overline{k} \to Y \times \mathrm{Spec}\,\overline{k}$ is smooth. In particular, if $X$ is a $k$-scheme $X$ is smooth over $k$ if and only if it is regular after changing the base to $\overline{k}$. With these in mind, the equivalence of the the two definitions is proved in

---
[1]Partially supported by NSF grant DMS-9401391

Theorem (1.8) [1, p. 130]. These definitions are equivalent to the one given in [2] also (see Theorem (10.2) [2, p. 269]).

An easy consequence of the second definition is that if $X$ is smooth over a field $k$, the set $X_{k^{\text{sep}}}$ is Zariski dense in $X$. Suppose that $k$ is algebraically closed, and both $X, Y$ are smooth over $k$. An important criterion of smoothness is that $f$ is smooth if and only if for all $x \in X$, the induced map between the tangent spaces $T_x X \to T_{f(x)} Y$ is surjective (see Corollary (5.5) [1, p. 149]). Also if $f : X \to Y$ is smooth and $X$ is smooth over $k$, $Y$ is smooth over $k$. One can see this in the following manner.

Let $f(x) = y$ and $h : X \to \text{Aff}^n \times Y$ ($\text{Aff}^n$ is the affine $n$-space) be an etale morphism in a neighborhood of $x$ such that the composition with the projection to $Y$ is $f$. Let $h(x) = (z, y)$, and $\text{pr}_1, \text{pr}_2$ be the projections from $\text{Aff}^n \times Y$ to $\text{Aff}^n, Y$. Then by Proposition (4.9) [1, p. 117],

$$\Omega^1_{X/k} \cong h^* \Omega^1_{(\text{Aff}^n \times Y)/k} \cong h^*(\text{pr}_1^* \Omega^1_{\text{Aff}^n/k} \oplus \text{pr}_2^* \Omega^1_{Y/k})$$

in a neighborhood of $x$. Since $f$ is surjective, it is easy to see that $\dim \Omega^1_{Y/k} \otimes k(y) = \dim Y$ for all $y \in Y$. Therefore, by the definition of smoothness in [2, p. 268], $Y$ is smooth over $k$.

Now we are ready to prove the theorem.

*Proof of Theorem.* We may assume that $k$ is separably closed. We first assume that $k$ is algebraically closed.

Let $F_1, \cdots, F_n$ be irreducible codimension one components of $V \setminus U$, $F = F_1 \cup \cdots \cup F_n$, and $W = V \setminus F$. Then $W$ is affine, $U \subset W$, and the codimension of $W \setminus U$ in $W$ is greater than one. Let $G = \text{Spec}\, R$, and $S = R^{G_w}$ (the ring of invariants). Then the geometric quotient $X = G/G_w$ exists, it is affine (in fact $X = \text{Spec}\, S$), and the map $f : G \to X$ is smooth (see Proposition 0.9, Amplification 1.3 [5, pp. 16,30]). The point here is that $S$ is finitely generated because of Haboush's theorem (see [5, p. 191]). By properties of smooth morphisms, $X$ is smooth over $k$. Because of the universal property of the quotient, the map $G \to U$ defined by $g \to gw$ factors through $X$.

Since $G, V$ are smooth over $k$, $\dim T_w U = \dim U = \dim V$ and $\dim T_e G = \dim G$ ($T_w U$, etc., are tangent spaces). Since $G \to X$ is a principal fiber bundle and the generic fiber is isomorphic to $G_w$, $\dim G = \dim X + \dim G_w = \dim U + \dim G_w$. Note that the map $X \to U$ is set-theoretically bijective, because $X$ is the geometric quotient. Therefore, $\dim G_w = \dim T_e G - \dim T_w U$. Since $G_w$ is smooth over $k$,

$$\dim G_w = \dim T_e G_w = \dim \text{Ker}(T_e G \to T_w U).$$

Therefore,
$$\dim \text{Ker}(T_e G \to T_w U) = \dim T_e G - \dim T_w U.$$

This proves that $T_e G \to T_w U$ is surjective.

Since the map $G \to U$ is $G$-equivariant, it induces surjective maps of tangent spaces. So that is the case for $X \to U$ also. Since both $X, U$ are smooth over $k$, this implies that $X \to U$ is smooth. Since the map $X \to U$ is of relative dimension zero, it is etale. Since this map is set-theoretically bijective, $X$ is isomorphic to $U$ and $U$ is affine. Since $W$ is smooth over $k$, it is normal in particular. Since the codimension of $W \setminus U$ in $W$ is greater than one, any regular function on $U$ extends to $W$. Since both $U, W$ are affine, $U = W$. This proves (1). The statement (2) is obvious if $k$ is algebraically closed.

We now consider arbitrary $k$ (but we may still assume $k$ is separably closed). Let $X = \operatorname{Spec} S$ again. Since $G_{w\overline{k}}$ is Zariski dense in $G_w$, it is easy to see that
$$(R \otimes \overline{k})^{G_{w\overline{k}}} \cong S \otimes \overline{k}.$$
This implies that $S$ is finitely generated over $k$ and $X \times \operatorname{Spec} \overline{k}$ is the geometric quotient $(G \times \operatorname{Spec} \overline{k})/(G_w \times \operatorname{Spec} \overline{k})$. The $k$-scheme $X$ obviously has the universal property of the quotient in the category of $k$-schemes. Therefore, the map $G \to U$ defined by $g \to gw$ still factors through $X$. By the previous step, the map $X \times \operatorname{Spec} \overline{k} \to U \times \operatorname{Spec} \overline{k}$ is an isomorphism. Therefore, the map $X \to U$ is an isomorphism also. So $U$ is affine and the statement (1) follows by the same argument as in the previous step.

Let $f$ be the map $G \to X$ again. Since $f$ is smooth, for any point $x \in U_k$, the fiber $f^{-1}(x)$ over $x$ is smooth over $k$. Since $k$ is separably closed, we may choose a $k$-rational point $g \in f^{-1}(x)$. Then $gw = x$, and this proves (2). $\square$

The reason why we proved the above theorem is the following. For simplicity, let's consider a pair $(G, V)$ as above such that $V$ is an irreducible representation. If $\operatorname{ch} k = 0$, a pair $(G, V)$ as above is called regular if there exists a relative invariant polynomial $\Delta(x)$ whose Hessian is not identically zero. This definition is equivalent to the existence of a generic point $w$ whose stabilizer is reductive. Let $V^{\operatorname{ss}} = \{x \in V \mid \Delta(x) \neq 0\}$. Because of the irreducibility, $\Delta(x)$ is essentially unique. So the regularity implies that $V^{\operatorname{ss}}_{\overline{k}}$ is a single $G_{\overline{k}}$-orbit. If $\operatorname{ch} k = 0$, $\overline{k} = k^{\operatorname{sep}}$, and Igusa showed in [3, p. 269] that
$$G_k \setminus V_k^{\operatorname{ss}} \cong \operatorname{Ker}(\operatorname{H}^1(k, G_w) \to \operatorname{H}^1(k, G))$$
($\operatorname{H}^1(k, G_w), \operatorname{H}^1(k, G)$ are Galois cohomology sets).

We considered the question of rational orbit decompositions for certain regular prehomogeneous vector spaces in [8], [4] in all the characteristics, and we had to worry about the separability. Because of the above theorem, all we have to do is to make sure that $G_w$ is smooth and reductive, and by exactly the same argument as in [3], we have the above interpretation of $G_k \setminus V_k^{\operatorname{ss}}$ in terms of Galois cohomology sets, and we no longer have to worry about the separability. This was the motivation of this paper. Of course, we still have to show that $G_w$ is reductive, but that is necessary anyway to determine the orbit space $G_k \setminus V_k^{\operatorname{ss}}$, and in many cases, it can be worked out by case by case analysis as in [6].

Finally, we point out that once we show that $G_w$ is reductive without knowing $Gw$ is open, the condition $\dim G_w = \dim G - \dim V$ implies $Gw$ is open as in the case of $\operatorname{ch} k = 0$. To see this, we may assume $k$ is algebraically closed. Then since $\overline{Gw}$ is integral, the set of non-singular points is a non-empty open set (see Theorem 5.3 [2, p. 33]). The geometric quotient $X = G/G_w$ still exists, it is smooth over $k$, and the map $X \to \overline{Gw}$ is bijective. So fibers of $G \to \overline{Gw}$ and fibers of $G \to X$ are set-theoretically the same. Then by Propositions 9.1, 9.7 [2, pp. 256, 257], $\dim \overline{Gw} = \dim G - \dim G_w = \dim V$. Since $Gw$ is a locally closed set (since it is constructible and irreducible) and $G$ acts transitively on it, it is open in $V$.

# References

[1] Altman, A., and S. Kleiman. *Introduction to Grothendieck duality*, volume 146 of *Lecture Notes in Mathematics*. Springer-Verlag, Berlin, Heidelberg, New York, 1970.

[2] Hartshonrn, R. *Algebraic geometry*, volume 52 of *Graduate Texts in Mathematics*. Springer-Verlag, Berlin, Heidelberg, New York, 1977.

[3] Igusa, J. On a certain class of prehomogeneous vector spaces. *J. of Pure and Applied Algebra*, 47:265–282, 1987.

[4] Kable, A.C., and A. Yukie. Prehomogeneous vector spaces and field extensions II. To appear in Invent. Math.

[5] Mumford, D., J. Fogarty, and F. Kirwan. *Geometric invariant theory*. Springer-Verlag, Berlin, Heidelberg, New York, 3rd edition, 1994.

[6] Sato, M., and T. Kimura. A classification of irreducible prehomogeneous vector spaces and their relative invariants. *Nagoya Math. J.*, 65:1–155, 1977.

[7] Sato, M., and T. Shintani. On zeta functions associated with prehomogeneous vector spaces. *Ann. of Math.*, 100:131–170, 1974.

[8] Wright, D.J., and A. Yukie. Prehomogeneous vector spaces and field extensions. *Invent. Math.*, 110:283–314, 1992.

Akihiko Yukie
Oklahoma State University
Mathematics Department
401 Math Science
Stillwater OK 74078-1058 USA
yukie@math.okstate.edu
http://www.math.okstate.edu/~yukie